\documentclass[10pt]{amsart}

\usepackage{amscd,amssymb}

\newtheorem{theorem}{Theorem}[section]
\newtheorem{proposition}[theorem]{Proposition}
\newtheorem{lemma}[theorem]{Lemma}
\newtheorem{corollary}[theorem]{Corollary}

\theoremstyle{definition}
\newtheorem{definition}[theorem]{Definition}
\newtheorem{example}[theorem]{Example}

\theoremstyle{remark}
\newtheorem{remark}[theorem]{Remark}
\newtheorem{question}[theorem]{Question}

\newcommand{\abs}[1]{\lvert#1\rvert}

\newenvironment{alphenum}%
{%
 \begin{enumerate}
}
{%
 \end{enumerate}%
}

\newenvironment{romenum}%
{%
 \begin{enumerate}
}
{%
 \end{enumerate}%
}

\newcommand{\N}{{\mathbb N}}
\newcommand{\Z}{{\mathbb Z}}
\newcommand{\Q}{{\mathbb Q}}
\newcommand{\R}{{\mathbb R}}
\newcommand{\C}{{\mathbb C}}

\renewcommand{\P}{{\mathbb P}}
\newcommand{\CP}{{\mathbb{CP}}}

\newcommand{\ie}{{\it i.e.}}
\newcommand{\cf}{{\it cf.\ }}
\newcommand{\eg}{{\it e.g.}}

\renewcommand{\a}{{\alpha}}
\renewcommand{\b}{{\beta}}

\newcommand{\mero}{{ \, -\negmedspace \to\, }}

\DeclareMathOperator{\vol}{{vol}}
\DeclareMathOperator{\Const}{{Const}}

\DeclareMathOperator{\length}{{length}}
\DeclareMathOperator{\area}{{area}}
\DeclareMathOperator{\sys}{{sys}}
\DeclareMathOperator{\stabsys}{{stabsys}}
\DeclareMathOperator{\pisys}{{\pi-sys}}
\DeclareMathOperator{\z2sys}{{\Z_{2}-sys}}

\DeclareMathOperator{\Image}{{Im}}
\DeclareMathOperator{\tors}{{torsion}}

\DeclareMathOperator{\id}{{id}}
\DeclareMathOperator{\Jac}{{Jac}}
\DeclareMathOperator{\Cyl}{{Cyl}}
\DeclareMathOperator{\Cone}{{Cone}}
\DeclareMathOperator{\SO}{{SO}}

\DeclareMathOperator{\interior}{{int}}
\begin{document}

\title[Volume and geometric inequalities] %
{Volume of Riemannian manifolds, geometric inequalities, and
homotopy theory}

\author[M.~Katz]{Mikhail G.~Katz}
\address{ UMR 9973, D\'{e}partement de Math\'{e}matiques,
Universit\'{e} de Nancy~1, B.P.~239, 54506 Vandoeuvre, France}
\curraddr{School of Mathematical Sciences, Sackler Faculty,
Tel Aviv University, Ramat Aviv, Tel Aviv 69978 Israel}
\email{katzm@math.tau.ac.il}

\author[A.~Suciu]{Alexander I.~Suciu}
\address{Department of Mathematics,
Northeastern University,
Boston, MA 02115, USA}
\email{alexsuciu@neu.edu}

\dedicatory{Dedicated to Mel Rothenberg on his $65^{\text th}$ birthday}

\thanks{To appear in the Rothenberg Festschrift
(M.~Farber, W.~L\"{u}ck, S.~Weinberger, eds.),
Contemporary Math., Amer. Math. Soc., Providence, RI.
Available at {\ttfamily http://xxx.lanl.gov/abs/ math.DG/9810172.}}


\begin{abstract}
We outline the current state of knowledge regarding geometric
inequalities of systolic type, and prove new results, including
systolic freedom in dimension~$4$.  Namely, every compact, orientable,
smooth $4$-manifold $X$ admits metrics of arbitrarily small volume
such that every orientable, immersed surface of smaller than unit area
is necessarily null-homologous in $X$.  In other words, orientable
$4$-manifolds are $2$-systolically free.  More generally, let
$m$ be a positive even integer, and let $n>m$.  Then all
manifolds of dimension at most $n$ are $m$-systolically free
(modulo torsion) if all $k$-skeleta, $m+1\leq k\le n$, of
the loop space $\Omega(S^{m+1})$ are $m$-systolically free.
\end{abstract}

\maketitle

\section{Introduction}
\label{sec:intro}
Our purpose here is both to present new results and to outline the
current state of knowledge regarding geometric inequalities of
systolic type.  One of the new results (proved in \S\ref{sec:sys4}) is
the following.

\begin{theorem}
\label{thm:sysfree4}
Every compact, orientable, smooth $4$-manifold $X$ admits metrics of
arbitrarily small volume such that every orientable, immersed surface
of smaller than unit area is necessarily null-homologous in $X$.
\end{theorem}

In other words, orientable $4$-manifolds are $2$-systolically free.  A
precise definition of systolic freedom will be given in
\S\ref{history}.  Roughly speaking, a $k$-systole of a Riemannian
manifold $(X,g)$ is the least volume of a non-trivial $k$-cycle in
$(X,g)$.  To say that $X$ is $k$-systolically free means that the
volume of $g$ imposes no upper bound upon the $k$-systole.  One can
view this result as a way of producing metrics for which the systole
and the mass in the middle dimension do not agree
(\cf Theorem~\ref{grom} and \cite{GMM}).

\begin{remark}
\label{FS}
Consider sufficiently small perturbations $g$ of the Fubini-Study
metric $g_0$ on $\CP^2$, scaled to the same volume as $g_0$.
M.~Gromov proved in \cite{G2} (Corollary, p.~309) that one can always
find embedded $2$-spheres, homologous to the projective line
$\CP^1\subset \CP^2$, of $g$-area no greater than that of $\CP^1$ for
the $g_0$ metric.  On the other hand, our Theorem~\ref{thm:sysfree4}
shows that there may be no such spheres if we go far away from the
Fubini-Study metric in the space of Riemannian metrics on $\CP^2$.  We
see that here local and global diverge.
\end{remark}

\begin{remark}
\label{freedman}
The conclusion of the theorem may well be as strong as possible, in
the sense that {\em non-orientable} surfaces of small area
(constrained by the volume of the ambient manifold) may always be
present.  M.~Freedman~\cite{Fr} raised the question of what happens
with the systoles defined using homology with $\Z_2$-coefficients.
Gromov~\cite{G5} conjectures that the systolic inequalities are true
in this case, \ie, our theorem above on systolic freedom has no
$\Z_2$-analogue.  As the example in Remark~\ref{mod2} illustrates,
$\Z_2$-cycles may have asymptotically smaller area growth than
$\Z$-cycles, supporting the conjecture that volume is a constraint
upon the $\Z_2$-systoles.  The exact boundaries between freedom and
constraint are yet to be described.
\end{remark}

In the $4$-dimensional case, we reduce the $\bmod\: 2$ problem to
$\CP^2$ and $S^2 \times S^2$, in the following sense.

\begin{theorem}
\label{nonorient}
All compact, orientable, smooth $4$-manifolds are $2$-systolically
free over $\Z_2$ if and only if $S^2 \times S^2$ and $\CP^2$ are
$2$-systolically free over $\Z_2$.
\end{theorem}

Our Theorem~\ref{thm:sysfree4} permits the following generalization of
the results obtained by the authors in collaboration with I.~Babenko
in \cite{BKS}.
\begin{theorem}
\label{middle}
Every compact, smooth manifold of dimension $n=2m\ge 4$ admits metrics
of arbitrarily small $n$-dimensional volume such that every
orientable, immersed middle-dimensional submanifold of smaller than
unit $m$-volume represents a torsion class in $m$-dimensional homology
with integer coefficients.
\end{theorem}

In other words, $2m$-dimensional manifolds are $m$-systolically free
(modulo torsion).  The term `systolic freedom (modulo torsion)' is
also explained in \S\ref{history}.  For manifolds of dimension larger
than $2m$, we prove the following.

\begin{theorem}
\label{thm:introbeyond}
Let $m\geq 2$ be an even integer, and let $n>m$.  Then all manifolds
of dimension at most $n$ are $m$-systolically free (modulo torsion) if
all $k$-skeleta, $m+1\leq k\le n$, of the based loop space
$\Omega(S^{m+1})$ are $m$-systolically free.
\end{theorem}

The paper is organized as follows.  In \S\ref{sec:iiss}, we define
systoles and give representative proofs of inequalities exhibiting
constraint imposed by the volume (of the ambient manifold) on them.
In \S\ref{history}, we explain the first example of systolic freedom
(\ie, the absence of such constraint), due to M.~Gromov, and give
related definitions and historical comments.
In \S\ref{subsec:beyond}, we prove Theorem~\ref{thm:introbeyond}.
In \S\ref{sec:carveup}, we introduce a `carving-up' technique, with the
goal of reducing the problem to the case when the middle Betti number
is at most $2$.
In \S\ref{sec:midsys}, we prove Theorem~\ref{middle}.
In \S\ref{sec:mero} we introduce a morphism, weaker than a continuous
map, suitable for pulling back systolic freedom.
In \S\ref{sec:sys4}, we prove Theorem~\ref{thm:sysfree4} and
Theorem~\ref{nonorient}.
In \S\ref{sec:twofreedom}, we prove a refinement of Theorem
\ref{thm:introbeyond} for the the $2$-systole.
In \S\ref{subsec:oddsys}, we establish the freedom of odd-dimensional
systoles.

The paper concludes with 5 appendices.
In Appendix~\ref{sec:compdim},  we describe the techniques used in
proving systolic freedom in the  context of a pair of
distinct complementary dimensions.
In Appendix~\ref{sec:freepull}, we prove a workhorse lemma,
which allows us to propagate systolic freedom.
In Appendix~\ref{sec:mv}, we prove a lemma (used in
paragraph~\ref{smm}) on surgeries along curves in $4$-manifolds.
In Appendix~\ref{sec:spin}, we use the technique of spinning to
identify explicitly the result of a useful surgery.
In Appendix~\ref{sec:wh}, we recall the necessary material on
Whitehead products and Hilton's theorem on homotopy groups of
bouquets of spheres.

\subsection*{Acknowledgments}
The first author is grateful to Misha Farber for the opportunity to
present an early version of this paper at the Topology conference at
Tel~Aviv University in June of 1998, and to Shmuel Weinberger for a
helpful discussion of the results herein.

\section{Isoperimetric inequalities and stable systoles}
\label{sec:iiss}
\subsection{Isoperimetric inequality}
\label{subsec:isoineq}
Every simple closed curve in the plane satisfies the inequality
\[
\frac{A}{\pi}\leq\left(\frac{L}{2\pi}\right)^2,
\]
where $L$ is the length of the curve and $A$ is the area of the region
it bounds.  This classical {\em isoperimetric inequality} is sharp
insofar as equality is attained only by a round circle.

\subsection{Loewner's theorem}
\label{subseq:loewner}
In the 1950's, C.~Loewner and P.~Pu proved the following two theorems.

Let $\R\P^2$ be the real projective plane endowed with an arbitrary
metric, \ie, an embedding in some $\R^n$.  Then
\[
\left(\frac{L}{\pi} \right)^2\leq \frac{A}{2\pi},
\]
where $A$ is its total area and $L$ is the length of its shortest
non-contractible loop.  This {\em isosystolic inequality} is also
sharp, to the extent that equality is attained only for a constant
curvature metric (quotient of a round sphere).

Similarly, every metric torus $T^2=\R^2/\Z^2$ satisfies the sharp
inequality
\[
L^2\leq \tfrac{2}{\sqrt{3}}A,
\]
where $L$ is the length of its shortest non-contractible loop, and $A$
is the area.

\subsection{What is a systole?}
In the 1970's, Marcel Berger initiated the study of a new Riemannian
invariant, which eventually came to be called the {\em systole}.

\begin{definition}  The (homology) $1$-systole
of a Riemannian manifold $(X,g)$ is the quantity
\[
\sys_1(g)=\inf_\ell\length(\ell),
\]
where the infimum is taken over all closed curves $\ell$ in $X$ which
are not homologous to zero.
\end{definition}

A similar homotopy invariant, $\pisys_1$, is obtained by minimizing
lengths of {\em non-contractible} curves.  These two invariants
obviously coincide for $\R\P^2$, $T^2$, and any manifold with abelian
fundamental group.

\subsection{Conformal representation and Cauchy-Schwartz}
We give a slightly modified version of M.~Gromov's proof of Loewner's
theorem for the $2$-torus, or rather the following slight
generalization.

\begin{theorem}
\label{thm:loewner}
There exists a pair of distinct closed geodesics on an arbitrary
metric $2$-torus $T$, of respective lengths $L_1$ and $L_2$, such that
\[
L_1 L_2 \leq \tfrac{2}{\sqrt{3}}A,
\]
where $A$ is the total area of the torus.  Equality is attained
precisely (up to scaling) for the equilateral torus $(1,\zeta)$, where
$\zeta$ is a primitive sixth root of unity.  Moreover, their homotopy
classes form a generating set for $\pi_1(T)= \Z\times\Z$.
\end{theorem}

\begin{proof}  A conformal representation $\phi: T_0 \rightarrow T$, where
$T_0$ is flat, may be chosen in such a way that $T$ and $T_0$ have the
same area.  Let $f$ be the conformal factor of $\phi$.

\begin{figure}
\begin{picture}(6,2.5)(0,0)
\multiput(1,0.5)(0.244,0.122){10}{\line(1,0){4}}
\multiput(1,0.5)(4,0){2}{\line(2,1){2.2}}
\put(0.5,0.5){\vector(0,1){1.122}}
\put(0.5,1.622){\vector(0,-1){1.122}}
\put(0.25,1.1){\makebox(0,0){$\sigma$}}
\put(3,0.2){\makebox(0,0){$\ell_0$}}
\put(6.8,1.1){\makebox(0,0){$\ell_s$}}
\put(6.55,1.13){\vector(-4,-1){1}}
\end{picture}
\caption{\textsf{Family of parallel geodesics on a flat torus}}
\label{fig:par}
\end{figure}

Let $\ell_0$ be any closed geodesic in $T_0$.  Let $\{\ell_s\}$ be the
family of geodesics parallel to $\ell_0$, see Figure~\ref{fig:par}.
Parametrize the family $\{\ell_s\}$ by a circle $S^1$ of length $\sigma$
so that $\sigma \ell_0 = \area(T_0)$.  Then $\area(T)=\int_{T_0}\Jac_\phi =
\int_{T_0}f^2$.  By Fubini's theorem, $\area(T)=\int_{S^1}ds
\int_{\ell_s} f^2 dt$.  By the Cauchy-Schwartz inequality,
\[
\area(T)\geq\int_{S^1}ds \frac{ (\int_{\ell_s} f dt)^2}{\ell_0 } =
\frac{1}{\ell_0}\int_{S^1} ds (\length\phi(\ell_s))^2.
\]
Hence there is an $s_0$ such that $\area(T)\geq
\frac{\sigma}{\ell_0} \length \phi(\ell_{s_0})^2$, \ie,
\[
\length\phi(\ell_{s_0}) \leq \ell_0.
\]

This reduces the proof to the flat case.  Given a lattice in $\C$, we
choose a shortest lattice vector $L_1$, as well as a shortest one
$L_2$ not proportional to $L_1$.  The inequality is now obvious from
the geometry of the standard fundamental domain in $\C$.
\end{proof}

\begin{theorem}[Gromov~\cite{G1}]
\label{grom}
Let $X^n$ be a compact, orientable, smooth manifold of dimension~$n$.
Let $\pisys_1(g)$ be the length of the shortest non-contractible loop
for the metric $g$ on $X$.  Then  the inequality
\begin{equation}
\label{eq:gromovineq}
(\pisys_1(g))^n\leq \Const\cdot\vol(g)
\end{equation}
holds for a positive constant $\Const$ independent of the metric, if
an only if the inclusion of $X$ in an Eilenberg-MacLane space
$K=K(\pi_1(X),1)$ retracts to the $(n-1)$-skeleton of $K$.  For
manifolds satisfying (\ref{eq:gromovineq}), moreover, the constant
$\Const=\Const_n$ depends only on $n=\dim X$.
\end{theorem}

The converse of this theorem was clarified by I.~Babenko~\cite{Ba}
(\cf Appendix~\ref{sec:freepull}, where similar techniques are used).
Note that the class of manifolds satisfying
inequality~\eqref{eq:gromovineq} includes aspherical manifolds, as
well as real projective spaces.

\subsection{Higher systolic invariants}
Let $X$ be a finite $n$-dimensional simplicial complex endowed with a
piecewise smooth metric $g$.  Let $k\leq n$.

\begin{definition}  Let $\alpha\in H_k(X;\Z)$.  Define
\[
\|\alpha\|=\inf_{M\in\alpha}\vol_k(M),
\]
where the infimum is taken over all piecewise smooth cycles $M$
representing the class $\alpha$.  Here the volume of a (smooth)
singular simplex is that of the pullback of the quadratic form $g$ to
the simplex, and we take absolute values of the coefficients to obtain
the volume of the cycle.  Also define a `stable norm' by
\[
\|\alpha \|_s=\lim_{q\rightarrow\infty}\frac{\|q\alpha\|}{q}.
\]
\end{definition}
\noindent
Clearly, we have $\|\alpha \|_s\leq \|\alpha \|$.

\begin{definition}
\label{defi3}
We define the following three systolic invariants for
the metric $g$ on $X$:
\begin{alphenum}
\item
The ordinary systole
\[
\sys_k(g)=\inf_{\alpha\not=0}\|\alpha\|.
\]

\item
The systole `modulo torsion'
\[
\sys^{\infty}_k(g)=\inf_{\alpha\not={\rm torsion} }\|\alpha\|.
\]

\item
The stable systole
\[
\stabsys_k(g)=\inf_{\alpha\not={\rm torsion} }\|\alpha\|_s.
\]
\end{alphenum}
\end{definition}
\noindent
Evidently, we have $\stabsys_k(g)\leq \sys^{\infty}_k(g)$.

\subsection{Calibration and stable systolic inequalities}
Higher systolic invariants can sometimes also be constrained by the
volume.  Such constraint is illustrated by the following theorem.

\begin{theorem}[Gromov~\cite{G4}] For every metric $g$ on $\CP^n$,
\[
\stabsys_2^n(g)\leq K_n \vol_{2n}(g),
\]
for a suitable constant $K_n$.
\end{theorem}

\begin{proof} The proof is a calibration argument.
Let $\alpha=[\CP^1]$ be the standard generator of
$H_2(\CP^n; \Z)=\Z$.  Let $\eta\in H^2(\CP^n; \Z)$
be the dual generator.  Let $\omega$ be any $2$-form
representing $\eta$.  Then
\[
1=\int_{\CP^n} \omega^{\wedge n}\leq K\|\omega \|^n \vol_{2n}(g).
\]
Hence $1\leq K\|\eta\|^n\vol_{2n}(g)=K\,\frac{\vol_{2n}(g)}{\|\alpha\|_s^n}$,
since the norms $\|\ \|$ in $H^2(X)$ and $\| \ \|_s$ in $H_2(X)$ are dual
by~\cite{Fe}, p.~394.

Thus $\|\alpha\|_s^n\leq K\vol_{2n}(g)$, and so
$\stabsys_2^n(g)\leq K\vol_{2n}(g)$, a stable systolic inequality.
\end{proof}

Another example (also due to Gromov) of a stable systolic inequality
is the following.  Let $X=S^p\times S^q$ endowed with an arbitrary metric
$g$.  Then
\[
\stabsys_p(g)\cdot \stabsys_q(g) \le \vol_{p+q}(g).
\]
J.~Hebda proved further results of this type in the 1980's.

Our objective is to show that in general, $k$-systoles are not
constrained by the volume, unlike the above theorems for $T^2$,
$\R\P^2$, and $\stabsys_2(\C\P^n)$.

\section{History and definitions}
\label{history}
\subsection{Question of freedom or constraint}
\label{subsec:freeornot}
M.~Berger asked in 1972 if the systolic invariants can be constrained
by the volume.  M.~Gromov reiterated this question in \cite{GLP}:

\begin{question}
\label{milliondollar}
Can one replace stable systoles by ordinary ones in the above
inequalities?
\end{question}

The question was asked again in an IHES preprint in 1992 (which
ultimately appeared as \cite{G3}).  Shortly afterwards, Gromov
described the first example of systolic freedom, \ie, the violation of
the systolic inequalities (see paragraph \ref{subsec:gromov}).  The
educated guess today is that, if one uses integer coefficients,
systolic freedom predominates as soon as one is dealing with a
$k$-systole for $k\geq 2$ (but see Remark~\ref{freedman}).

\subsection{Gromov's example}
\label{subsec:gromov}
Gromov described metrics on $S^1\times S^3$ which provided an
unexpected negative answer to the question \ref{milliondollar}.  He
stated it in global Riemannian terms.  We provide a reformulation in
terms of differential forms, which lends itself easier to
generalization.

Let $b=J(dr)|_{TS^3}$ be the contact $1$-form on $S^3$, where
$S^3\subset\C^2$ is the unit sphere, $J$ is rotation by $\pi/2$
furnished by the complex structure, $dr$ is the $1$-form on
$\C^2\setminus\{0\}$, defined by the radial coordinate $r$.

Complete $b$ to an orthonormal basis $\{b,b',b''\}$ of $T^*S^3$, which
is canonically identified with $TS^3$ by means of the standard unit
sphere metric.

Let $S^1$ be the unit circle parametrized by $e^{iz}$, where $z$ is
real.  The standard $1$-form $dz$ is the arc-length.  Gromov's
sequence of left-invariant metrics, $\{g_j\}_{j=1}^{\infty}$, on
$S^1\times S^3$ is defined by the following quadratic forms:
\[
g_j=(dz-jb)^2 +b^2 + (1+j^2)({b'}^2+{b''}^2).
\]
The essential ingredient here is the matrix
\[
\begin{bmatrix}
1&-j\\ -j & 1+j^2
\end{bmatrix}
\begin{matrix}
\leftarrow  z \\
\leftarrow  b
\end{matrix}
\]
The coefficient $1+j^2$ ultimately determines the $3$-systole, the
coefficient $1$ determines the $1$-systole, and the determinant of the
matrix determines the volume of $g_j$.  Then we obtain
\[
\frac{\vol(g_j)}{\sys_1(g_j)\sys_3(g_j)}\xrightarrow[j \rightarrow
\infty]{} 0.
\]

See Appendix~\ref{sec:compdim} for a `local' version of this
example, suitable for generalization.

\begin{definition}
\label{def:mfree}
A finite $n$-dimensional CW-complex $X$ is {\em $m$-systolically free}
if
\[
\inf_{g} \frac{\vol_{n}(g)^\frac{m}{n}}{\sys_m (g)}=0,
\]
where the infimum is taken over all piecewise smooth metrics $g$ on a
simplicial complex $X'$ in the homotopy type of $X$.  
\end{definition}

\begin{remark}
\label{rem:sysfree}
That every finite CW-complex $X$ is homotopy equivalent to a finite simplicial complex 
$X'$ is well-known, see~\cite{LW}.  That the definition does not depend on the 
choice of $X'$ follows from Corollary~\ref{invariant} below.  The systolic freedom 
of $X$ is equivalent to the existence of an {\em $m$-free sequence of metrics}, 
$\{g_j\}$, on $X'$:
\[
\sys_m^n(g_j)\ge j\, \vol_{n}^m(g_j).
\]
\end{remark}

\begin{definition}
\label{def:sysmodtors}
An $n$-dimensional CW-complex is {\em $m$-systolically free (modulo torsion)} if
\[
\inf_{g} \frac{\vol_{n}(g)^\frac{m}{n}}{\sys^{\infty}_m (g)}=0.
\]
\end{definition}

Theorem~\ref{nonorient} in the Introduction uses the following
terminology in the context of coefficients modulo~$2$.

\begin{definition}
\label{def:z2free}
Let $X$ be a compact, smooth manifold of dimension~$n$.  Consider $m$-cycles with 
$\Z_2$-coefficients, \eg, maps of manifolds (orientable or not) into $X$,
which represent non-trivial classes in $H_m(X;\Z_2)$, and calculate
the $m$-volume of the pullback metric.  Define $\z2sys_{m}(X)$ to be
the infimum of all such volumes.  We say that $X$ is {\em
$m$-systolically free over $\Z_2$} if
\[
\inf_{g}\frac{\vol_{n}(g)^\frac{m}{n}}{\z2sys_{m}(g)}= 0.
\]
\end{definition}

\section{Reduction to loop space}
\label{subsec:beyond}
For even $m$, we reduce the systolic problem to spaces each of whose
Betti numbers is at most $1$.  More specifically, we show that the
systolic $m$-freedom of $n$-manifolds (modulo torsion) reduces to that
of certain skeleta of loop spaces of spheres.  Thus, whether or not
there exists an analogue of Gromov's Theorem~\ref{grom} for the higher
systoles, depends on the loop space in question.

Our proof of Theorem~\ref{thm:introbeyond} in the Introduction is
modeled on the proof of Lemma \ref{lem:prebouquet} below.

\begin{lemma}
\label{cohotrick}
Let $m\geq 1$ be an integer.  Let $X$ be a finite CW-complex.  Let
$b=b_m(X)$.  Let $K=K(\Z^{b},m)$ be an Eilenberg-MacLane space with
skeleta $K^{(m+1)}=K^{(m)}=\vee^b S^m$.  There is a map $X\to K$,
inducing an isomorphism $H_m(X,\Z)/{\rm torsion}\to H_m(K,\Z)$,
and whose restriction to $X^{(m+1)}$ has image in $\vee^b S^m$.
\end{lemma}

\begin{proof}
A basis for $H^m(X;\Z)/{\tors}$ defines a map $X\to K$, inducing the
required isomorphism.  The cellular approximation theorem yields a map
$X^{(m+1)}\to \vee^b S^m$.
\end{proof}

\begin{lemma}
\label{lem:prebouquet}
Any $(2m-1)$-dimensional CW-complex $X$ admits a map to a suitable
bouquet of $m$-spheres which induces an isomorphism in rational
homology of dimension~$m$.
\end{lemma}

\begin{proof}
Consider the bouquet $\vee^b S^m$, where $b=b_m(X)$.  According to
B.~Eckmann~\cite{E}, a map $\phi_q:S^m\to S^m$ of degree $q$ induces
multiplication by $q$ in the stable groups $\pi_j(S^m)$, for $m+1\leq
j\leq 2m-2$.  Thus, if $q$ is a multiple of the order of $\pi_j(S^m)$,
the self-map $\vee^b \phi_q:\vee^b S^m\to \vee^b S^m$ induces the zero
homomorphism in $\pi_j$ by Hilton's theorem \ref{subsec:hm}.  Hence a
map $f_j: X^{(j)}\to \vee^b S^m$, followed by $\vee^b \phi_q$, extends
to $f_{j+1}:X^{(j+1)}\to \vee^b S^m$.  The lemma now follows by
induction on $j$, based on Lemma~\ref{cohotrick}.
\end{proof}

\begin{corollary}
\label{pathological}
Let $m$ and $n$ be integers such that $2\leq m < n < 2m$.  Then every
$n$-dimensional manifold is $m$-systolically free (modulo torsion).
\end{corollary}

\begin{proof}
The bouquet of $m$-spheres, viewed as an $n$-dimensional complex, is
obviously $m$-free (the numerator in Definition \ref{def:mfree}
vanishes).  We apply the pull-back Lemma~\ref{lem:longcyl2}.
\end{proof}

\begin{lemma}
\label{lem:bouquet}
The $(2m-1)$-skeleton of a $(2m)$-manifold $X$ admits a map to a
suitable bouquet of $m$-spheres, inducing an isomorphism in
$m$-dimensional rational homology, and sending the attaching maps of
the $(2m)$-cells to a sum of Whitehead products.
\end{lemma}

\begin{proof}
We continue with the notation of the proof of
Lemma~\ref{lem:prebouquet}.  We thus have a map
$f_{2m-1}:X^{(2m-1)}\to \vee^b S^m$.  Let $e=[S^m]$ be the fundamental
class of the sphere.  Recall that the Whitehead product $[e,e]\in
\pi_{2m-1}(S^m)$ generates precisely the kernel of the suspension
homomorphism.  Suspension commutes with $\phi_q$.  Hence, if $q$ is a
multiple of the order of the stable group $\pi_{2m}(S^{m+1})$, then
$(\phi_q)_*(\pi_{2m-1}(S^m))$ is contained in the subgroup generated
by Whitehead products.  Hence $(\vee^b \phi_q)\circ f_{2m-1}$ is the
desired map.
\end{proof}

\begin{proof}[Proof of Theorem~\ref{thm:introbeyond}]
Consider the based loop space $L = \Omega(S^{m+1})$.  As shown by
J.-P.~Serre~\cite{Se}, all homotopy groups of $L$ are finite, except
$\pi_m$ (here we need $m$ even).  Let $q$ be the product of the orders
of the homotopy groups of $L$ from $m+1$ up to dimension $n$.  Let
$\psi:L\to L$ be the map that sends a loop $\omega$ to $\omega^q$.
Since the usual group structure on $\pi_i(L)$ coincides with the one
coming from the multiplication of loops in $L$, the map $\psi$ induces
multiplication by $q$ in each homotopy group of $L$.  Thus, $\psi = 0$
on $\pi_i(L)$, for $m+1\le i\le n$.

Now let $X$ be an $n$-manifold, and let $b=b_m(X)$.  Let $L^{\times b}$ 
be the Cartesian product of $b$ copies of $L$. With respect to a
suitable cell structure, the $m$-skeleton of $L^{\times b}$ is a bouquet
$\vee^b S^m$ of $b$ copies of $S^m$.  By Lemma \ref{cohotrick}, the
$(m+1)$-skeleton of $X$ admits a map $f_{m+1}$ to the $m$-skeleton of
$L^{\times b}$ which induces rational isomorphism in $m$-dimensional homology.

We now proceed as in the proof of Lemma~\ref{lem:prebouquet}, with the
bouquet of spheres replaced by $L^{\times b}$.  By induction, the map
$\psi\circ f_{k}$ extends to
\[
f_{k+1} : X^{(k+1)}\to L^{\times b}
\]
for all $k$ from $m+1$ to $n$.  In this fashion, we obtain a map
$f:X\to L^{\times b}$ such that $f_*:H_m(X;\Q)\to H_m(L^{\times b};\Q)$ is an
isomorphism.  By the pull-back Lemma~\ref{lem:longcyl2}, the freedom
of $X$ follows from that of $L^{\times b}$.  It remains to reduce 
the freedom of $L^{\times b}$ to that of $L$ itself.

A well-known result of I.~James (see~\cite{Wh}), states that $L$ has,
up to homotopy, a cell decomposition
\[
L= S^m \cup e^{2m} \cup e^{3m} \cup \cdots ,
\]
with precisely one cell in each dimension divisible by $m$.
Therefore, we may proceed as in \S\ref{sec:carveup}, and carve up the
$n$-skeleton of $L^{\times b}$, reducing the problem to the freedom of the
closures of $n$-dimensional cells of $L$.  Each such closure is a
product of skeleta of $L$.  Therefore, Theorem~\ref{thm:introbeyond}
results from the following lemma.
\end{proof}

\begin{lemma}
\label{further}
Let $m\geq 2$.  Then the product of an $(m-1)$-connected,
$m$-systolically free complex with $S^m$ is again $m$-free.
Furthermore, the product of $(m-1)$-connected, $m$-systolically free
complexes is again $m$-free.
\end{lemma}

\begin{proof}
Given $m$-free families of metrics on each factor, we scale them to
unit $m$-systole.  We also scale the $m$-sphere to unit $m$-volume.
The product metrics then form $m$-free families.
\end{proof}

\section{Carving-up procedure}
\label{sec:carveup}
The procedure in question is helpful in breaking up the problem of
systolic freedom of complicated spaces into simpler components.  We
introduce it here as a way of streamlining the arguments of
\cite{BKS}, where the authors, in collaboration with I.~Babenko,
proved the middle-dimensional systolic freedom of even-dimensional
manifolds of dimension $n=2m\geq 6$ with free $H_m(X,\Z)$.

Let $K$ be an $n$-dimensional CW-complex with no cells in dimension
$n-1$.  The carving-up procedure consists of inserting a cylinder
between each $n$-dimensional cell and its boundary.  The new complex
$K'$ has the following three properties.

\begin{enumerate}
\item $K'$ has the same homotopy type as $K$.
\item $K'$ has the same number of top-dimensional cells $e^n_i$ as $K$.
\item The closures $\bar{e}^n_i$ of all top-dimensional cells in $K'$
are {\em disjoint}.
\end{enumerate}
Assume for the sake of simplicity that the image of each attaching map
\[
f_i: \partial D^n \rightarrow K^{(n-2)}
\]
of $e_i$ is a subcomplex (\ie, there are no partly covered cells in
the $(n-2)$-skeleton).  Let
\[
\partial^i=\Image(f_i)\subset K^{(n-2)}
\]
be the boundary of the $i^{\text{th}}$ top-dimensional cell.

The procedure of inserting cylinders in $K$ is formalized as follows:
\[
K'= K^{(n-2)} \cup_{g^0_i} (\partial^i \times [0,1])
\cup_{g^1_i} \bar e^n_i.
\]
Here $K^{(n-2)}$ is the $(n-2)$-skeleton, and the extremities of the
$i^{\text{th}}$ cylinder are attached along the inclusion maps $g^0_i:
\partial^i \times \{0\} \rightarrow K^{(n-2)}$, and $g^1_i: \partial^i
\times \{1\} \rightarrow \bar{e}^n_i$.  The repeated upper and lower
index $i$ indicates that the procedure is repeated for each
top-dimensional cell, as in Einstein notation.  Thus, we insert as
many cylinders as there are top-dimensional cells, to obtain $K'$.

\begin{lemma}
\label{lem:closures}
Let $K$ be an $n$-dimensional CW-complex with no cells in dimension
$n-1$.  Then the systolic $q$-freedom of $K$ reduces to that of the
closures of its top-dimensional cells.
\end{lemma}

\begin{proof}
We replace $K$ by $K'$ as above.  To make sure that areas and volumes
are well-defined, we replace each $\partial^i$ by a simplicial complex
of the same dimension and in the same homotopy type.  Then the
positive codimension condition is satisfied.  We also replace each
cell closure by a simplicial complex in the same homotopy type so as
to extend the simplicial structure on $\partial^i$.  We make a similar
replacement for $K^{(n-2)}$, obtaining a simplicial complex $K''$.

We apply to $K''$ the long cylinder argument in the proof of
Lemma~\ref{lem:longcyl1} and Lemma~\ref{lem:longcyl4}.  The key tools
here are the isoperimetric inequality and the coarea inequality,
applied to $\partial^i\times [0,1]$, where the metric on the interval
is chosen sufficiently long, to minimize interaction between opposite
extremes of the cylinder.
\end{proof}

\section{Middle-dimensional freedom}
\label{sec:midsys}
We now establish systolic freedom (modulo torsion) in the middle
dimension, as stated in Theorem~\ref{middle}.  In fact, we show that
the theorem is valid in the more general case of triangulable spaces
with piecewise smooth metrics, for which of course areas and volumes
can still be defined.

In contrast with Loewner's theorem \ref{thm:loewner}, the
middle-dimen\-sional systole ($m$-systole) is not constrained by the
volume when $m\geq 2$.  Thus Loewner's theorem has no
higher-dimensional analogue if one works with orientable submanifolds.
On the other hand, there might be such an analogue if one allows
non-orientable submanifolds (see Remark~\ref{freedman}).

A proof in the case $m\ge 3$ appeared in \cite{BKS}.

\subsection{Idea of proof}
\label{subsec:idea}
The idea of the proof of Theorem~\ref{middle} is to reduce the problem
to a local version of Gromov's example, described in
\S\ref{subsec:gromov}, by using pullback arguments of
Appendix~\ref{sec:freepull} as follows.  High-degree self-maps of
$S^m$ combined with Hilton's theorem \ref{subsec:hm} allow us to map
$X$ to a kind of a first-order approximation to the $2m$-skeleton of
the Eilenberg-MacLane space $K(\Z^b, m)$.  This approximation has the
same rational homology and contains the minimal number of
$2m$-dimensional cells.  The carving-up procedure of
\S\ref{sec:carveup} reduces the problem to the freedom of the closures
of these cells, each with middle Betti number at most $2$.  An
additional pull-back reduces the problem to the freedom of $S^m \times
S^m$.  Finally, the latter is reduced to the $(1,m)$-freedom of
$S^1\times S^m$ (\cf~Appendix~\ref{sec:compdim}).

\subsection{Proof of Theorem~\ref{middle}}  Let us restate
the theorem in a more convenient (and slightly more general) setting.

\begin{theorem}
\label{thm:midfreeCW}
Let $m\geq 2$.  Let $X$ be a finite, triangulable CW-complex of
dimension $2m$.  Then $X$ is $m$-systolically free (modulo torsion).
\end{theorem}

\begin{proof}
Let $b=b_m(X)$.  Let $P=S^m\cup_{[e,e]}D^{2m}$, where $e$ is the
fundamental class of $S^m$.  Consider the Cartesian product
$Q=P^{\times b}$ of $b$ copies of $P$.  Its $m$-skeleton
$Q^{(m)}=\vee^{b} S^m$ is a bouquet of $b$ copies of $S^m$.  Let
$f:X^{(2m-1)}\to Q^{(m)}$ be the map given by Lemma~\ref{lem:bouquet}.
Here the image of $\pi_{2m-1}(X^{(2m-1)})$ is contained in the
subgroup generated by Whitehead products.  Then $f$, followed by the
inclusion $Q^{(m)}\hookrightarrow Q^{(2m)}$, extends across all of $X$
by Hilton's theorem \ref{subsec:hm}.  Thus the freedom of $X$ reduces
to that of the $(2m)$-skeleton of $P^{\times b}$.

We now apply Lemma~\ref{lem:closures} to obtain a further reduction to
only two cases: $S^m\times S^m$ and $P$ itself.
The CW-complex $P$ is homotopy equivalent to the regular CW-complex
$W=S^m \times S^m \cup_{a+b} D^{m+1}$ (product of spheres with a disk
attached along the diagonal).  Thus the systolic $m$-freedom of $P$
reduces to that of $S^m \times S^m$ by Lemmas \ref{lem:longcyl1}
and \ref{lem:longcyl2}.

Finally, the middle-dimensional freedom of $S^m\times S^m$ follows
from Lemmas \ref{lem82} and \ref{lem:frees2t2} below.
\end{proof}

\begin{lemma}
\label{lem82}
Let $m\geq 2$.  Then the manifold $S^m\times S^m$ admits a map to a
CW-complex obtained from $S^1\times S^{m-1}\times S^m$ by attaching cells
of dimension at most $2m-1$, which induces a monomorphism in
$m$-dimensional homology.
\end{lemma}

\begin{proof}
The manifolds $S^m$ and $S^1\times S^{m-1}$ are related by surgery,
and thus $S^m \times S^m$ admits a map to $S^1 \times S^{m-1}\times
S^m\cup D^2\times * \times S^m$.

For $m=2$, this does not define a meromorphic map, since we have added
top dimensional cells.  In this case, we proceed as in
Corollary~\ref{cor:s2s2mero}.
\end{proof}

\begin{lemma}
\label{lem:frees2t2}
The manifold $S^1 \times S^{m-1}\times S^m$ is $m$-systolically free
if $m\geq 2$.
\end{lemma}

\begin{proof}
The manifold $S^1 \times S^m$ is $(1,m)$-free by
Theorem~\ref{thm:compfree}.  Let $g_j$ be such a free sequence
of metrics on $S^1 \times S^m$.  Taking the product with a sphere
$S^{m-1}$ of volume equal to $\frac{\sys_m(g_j)}{\sys_1(g_j)}$, we
obtain an $m$-free sequence of metrics on $S^1 \times S^{m-1}\times
S^m$.
\end{proof}

\subsection{Spin manifolds}
\label{subsec:spinmfd}
Our Theorem~\ref{thm:sysfree4} improves the general middle-dimen\-sional
result in the case $m=2$, to the extent that it removes the `modulo
torsion' clause.  A similar improvement exists for $m=3$ and $m=4$.

\begin{proposition}
Spin manifolds of dimension $6$ and $8$ are systolically free in
middle dimension.
\end{proposition}

\begin{proof}
We reduce the problem to the simply-connected case as in
paragraph \ref{subsec:orientfree}.  Furthermore, the spin condition
$w_2=0$ ensures that all embedded $2$-spheres have trivial normal bundles.
Let $C$ be the union of a disjoint family of embedded $2$-spheres
representing a set of generators for $2$-dimensional homology, and
perform surgery along $C$.  An analogue of Lemma~\ref{prop:merosurg2}
reduces the problem to the $2$-connected case.  A $2$-connected $6$-manifold
has torsion-free $3$-dimensional homology, and the proposition is
established for a $6$-manifold.

For $8$-manifolds, we continue by choosing a disjoint family of
embedded $3$-spheres (whose normal bundles are automatically trivial)
which represent generators for $3$-dimensional homology, and argue
as before, reducing the problem to the $3$-connected case.
\end{proof}

\section{Meromorphic maps}
\label{sec:mero}
Here we develop a convenient language for establishing the
$2$-systolic freedom of orientable $4$-manifolds.  We define a
morphism, weaker than a continuous map, suitable for pulling back
systolic freedom.  Our technique is introduced most provocatively by
means of the following question.

\begin{question}
What do meromorphic maps and surgeries have in common?
\end{question}

The answer is, roughly speaking, as follows (see
Example~\ref{exm:blowup} and Lemma \ref{prop:merosurg2}): both of
them carry an underlying structure, crystallized in the concept of a
`topological meromorphic map' below.

\subsection{`Topological meromorphic maps'}
Here we attempt to define a morphism, weaker than a continuous map,
suitable for pulling back systolic freedom.  Such morphisms could also
be called `topological blow-up maps' or `topological rational maps'.

\begin{definition}
\label{def:mero}
Let $X^n$ and $Y^n$ be manifolds of dimension $n=2m$.  A ``meromorphic
map", denoted by a broken arrow below,
\[
X\mero Y
\]
is a continuous map $f:X\to W$, such that:

\begin{alphenum}
\item\label{mero1}
The space $W$ has the homotopy type of a CW-complex obtained from
$Y$ by attaching cells of strictly smaller than the top dimension:
$W\simeq Y\cup \bigcup_{i} e^{\leq n-1}_{i}$.
\item\label{mero2}
The map $f$ induces a monomorphism in the middle dimension:
$\ker(f_*: H_m(X)\rightarrow H_m(W)) = 0.$
\end{alphenum}
\end{definition}

\begin{remark}
\label{epi}
We do not require the inclusion $Y\to W$ to induce an epimorphism in
middle-dimensional homology.
\end{remark}

\begin{example}
\label{exm:blowup}
Let $X$ be a complex manifold of (complex) dimension~$m$, let
$\widehat{X} \to X$ be the blow-up at a point $p\in X$, and let $\phi:
X\mero \widehat{X}$ be the classical meromorphic map (undefined at
$p$).  For $m=2$, set $W=\widehat{X}\cup D^{3}$, where the $3$-cell is
attached along the exceptional curve $S^2$.  Then $\phi$ can be
modified in a neighborhood $D^4$ of $p$ so as to lift a continuous map
$f: X\to W$ (a homotopy equivalence).  Here we take the cone $D^4\to
D^3$ of the Hopf fibration $S^3\to S^2$.  More precisely, the Hopf
fibration extends to a continuous map whose restriction to the
complement of $D^4$ in $X$ is a diffeomorphism onto the complement of
$S^2$ in $\widehat{X}$, while $D^4$ is mapped to $D^3$ as described.
For general $m$ we have a continuous map $X\to \widehat{X}\cup
\Cone(\CP^{m-1})$.
\end{example}

\subsection{Surgery and meromorphic maps}
\label{smm}
We now show that, under certain homological conditions, surgeries
along curves yield meromorphic maps between $4$-manifolds.

\begin{lemma}
\label{prop:merosurg2}
Let $X$ be a closed, orientable, smooth $4$-manifold.  Let $C\subset
X$ be a union of smoothly embedded, disjoint closed curves.  Let $Y$
be the result of surgery along $C$.  Then $X$ admits a meromorphic map
to $Y$.  The conclusion holds equally well if one uses coefficients
modulo~$2$ in the definition of homology and meromorphic maps.
\end{lemma}

\begin{proof}
Since $X$ is orientable, the normal disk bundle of $C$ is
trivializable.  Over each component of $C$, there are two
possible trivializations (or, {\em framings}), corresponding
to $\pi_1(\SO(3))=\Z_2$.  The framing (which is part of the
surgery data) identifies the normal disk bundle with $C \times D^3$.
By definition,
\[
Y=(X \setminus C \times \interior D^3) \cup {\textstyle \bigcup_{i}}
(D^2_{i} \times S^2),
\]
where $\bigcup_{i} D^2_{i}$ is a disjoint union of $2$-disks
(one for each  connected component $C_i$ of $C$).
Let $Z=X \times I \cup \bigcup_{i} D^2_{i}\times D^3$
be the result of attaching a handle of index~$2$ to $X$
along each component of $C$, according to the given framing.
Then $Z$ is the cobordism between $X$ and $Y$ determined by the
surgery, see J. Milnor~\cite{Mi}.

Now let $W=X\cup \bigcup_{i} D^2_{i}$ be the mapping cone of the
inclusion $C\subset X$.  In other words, $W$ is the CW-complex
obtained from $X$ by attaching the cores of the handles.

Clearly, the inclusion of $X$ in $W$ induces a monomorphism in
$2$-dimensional homology.  Moreover, $W$ is a deformation-retract of
$Z$, see~\cite{Mi}.  Figure~\ref{fig:cobordism} depicts the various
spaces introduced so far.

\begin{figure}
\begin{picture}(6,4)(0,0.3)
\put(1,1){\framebox(4,2.5){}}
\put(3,1){\oval(2,2)[t]}
\put(3,1){\oval(1.5,1.5)[t]}
\put(3,1){\oval(1,1)[t]}
\multiput(2.25,1)(1.5,0){2}{\circle*{0.15}}
\put(4.5,0.75){\makebox(0,0){$X$}}
\put(4.5,3.75){\makebox(0,0){$Y$}}
\put(0.7,2.5){\makebox(0,0){$Z$}}
\put(3,0.4){\makebox(0,0){$C$}}
\put(3.1,0.65){\vector(2,1){0.5}}
\put(2.9,0.65){\vector(-2,1){0.5}}
\put(4,2.25){\makebox(0,0){$D^2$}}
\put(3.6,2.05){\vector(-2,-1){0.5}}
\put(1.7,1.7){\makebox(0,0){$W$}}
\put(1.85,1.6){\vector(2,-1){0.4}}
\put(1.65,1.5){\vector(0,-1){0.44}}
\end{picture}
\caption{\textsf{Surgery and cobordism}}
\label{fig:cobordism}
\end{figure}

We have an inclusion $i:X\subset Z$, as well as a homotopy equivalence
$Z\simeq Y\cup e^3$.  This gives a map $X\to Y\cup e^3$, which
satisfies condition~(\ref{mero2}) in Definition~\ref{def:mero} of
injectivity in $2$-dimensional homology, since $Z$ is homotopy
equivalent to $W$.
\end{proof}

\begin{corollary}
\label{simplify}
An orientable 4-manifold admits a meromorphic map to a simply
connected one.
\end{corollary}

\begin{proof}
Let $C\subset X$ be the union of a family of disjoint embedded closed
curves representing a set of generators for the fundamental group.
Surgery on $X$ along $C$ produces a simply-connected manifold $Y$, and
we apply  Lemma \ref{prop:merosurg2}.
\end{proof}

Of particular interest to us (\cf Corollary~\ref{cor:s2s2mero}) is a
`dual' version of Lemma~\ref{prop:merosurg2} (the meromorphic map goes
the other way).

\begin{proposition}
\label{lem:merosurg}
Let $X$ be a closed, orientable, smooth $4$-manifold.  Let
$Y$ be the result of surgery along a union $C$ of smoothly
embedded, disjoint closed curves in $X$.  Assume that the connected
components of $C$ define a linearly independent set in $H_1(X;\Q)$.
Then $Y$ admits a meromorphic map to $X$.
\end{proposition}

\begin{proof}
We continue with the notation of Lemma~\ref{prop:merosurg2}.
Let $j: Y\to Z$ be the inclusion and $g: Z\to W$, the
homotopy equivalence.   Then $g_{*}\circ j_{*}:H_2(Y)\to H_2(W)$
is a monomorphism.  This is established by induction on
the number of connected components of $C$,
using Lemma~\ref{lem:circsurg} in Appendix~\ref{sec:mv}.

More precisely, let $\psi$ be the isomorphism provided by
Lemma~\ref{lem:circsurg}.  Since $W$ is obtained from $X$ by attaching
$2$-cells, the inclusion $i:X\to Z$ induces a monomorphism
$i_*:H_2(X)\to H_2(Z)$.  Since $j_*\circ \psi = i_*$, it follows that
$j_*$ is injective, and so $g_{*}\circ j_{*}$ is, too.

Hence, $g\circ j:Y\to W$ satisfies condition~(\ref{mero2})
in Definition~\ref{def:mero}.  We thus have the required
meromorphic map $Y\mero X$.
\end{proof}

\begin{remark}
\label{rem:linindep}
Without the independence hypothesis, $Y$ may have a larger second
Betti number than either $X$ or $W$.  For example, if $X=S^4$ and
$C=S^1$, then $W\simeq S^2 \vee S^4$ and $Y=S^2\times S^2$ or
$Y=S^2\widetilde{\times} S^2\cong \CP^2\# \overline{\CP}^2$, depending
on the parity of the framing.  Thus $b_2(X)=0$, $b_2(W)=1$, and
$b_2(Y)=2$.
\end{remark}

\begin{corollary}
\label{cor:s2s2mero}
There exists a meromorphic map $S^2 \times S^2 \mero S^2 \times T^2$.
\end{corollary}

\begin{proof}
Surgery on $T^2\times S^2$ along a pair of generators of $\pi_1(T^2
\times S^2)=\Z\times \Z $ yields $S^2 \times S^2$.  An explicit
verification of this fact is given in Corollary~\ref{cor:spin}.
Alternatively, one readily sees that the surgery does not change the
intersection form, and the resulting manifold is simply-connected.
Hence it is certainly homotopy equivalent to $S^2\times S^2$.  Either
way, Proposition~\ref{lem:merosurg} applies, completing the proof.
\end{proof}

\subsection{Pullback lemma for free metrics}
\label{subsec:pullbacks}
The following proposition allows us to propagate the phenomenon of
freedom once we exhibit it for products of spheres.

\begin{proposition}
\label{prop:xyfree}
Let $f: X \mero Y$ be a meromorphic map.  If $Y$ is systolically free,
then so is $X$.
\end{proposition}

\begin{corollary}  
\label{invariant}
Systolic freedom is a homotopy invariant.
\end{corollary}

\begin{proof}  Any map defining a homotopy equivalence is obviously a
``meromorphic map" (here $W=Y$, \ie, the set of attached cells is
empty).
\end{proof}

Proposition~\ref{prop:xyfree} follows from the following
lemma (\cf Lemma 6.1 of \cite{BaK}).

\begin{lemma}
\label{lem:freep}
Let $X$ and $Y$ be triangulable CW-complexes of dimension $n$. Suppose
there is a map $\phi: X \to W$ inducing a monomorphism on $H_q$, where
$W=Y\cup_{f} e^{k}$ is obtained from $Y$ by attaching a single cell of
dimension $k\leq n-1$.  Then the $q$-systolic freedom of $Y$ implies
that of $X$.  Moreover, if $X$ is a smooth manifold, the metrics
exhibiting freedom may be chosen to be smooth.
\end{lemma}

The lemma is proved in Appendix~\ref{sec:freepull}.

\section{Systolic freedom in dimension~$4$}
\label{sec:sys4}
In this section, we show that orientable closed $4$-manifolds are
systolically free in the middle dimension.

In the simply-connected case, the first author, in collaboration with
I.~Babenko \cite{BaK}, already reduced the problem to the case of
$S^2\times S^2$.  Here we notice that $S^2\times S^2$ admits a
meromorphic map to $S^2\times T^2$ (see Corollary~\ref{cor:s2s2mero}),
and so its middle-dimensional freedom results from that of
$S^2\times T^2$ (see Lemma~\ref{lem:frees2t2}).

\subsection{Systolic freedom modulo torsion}
\label{subsec:modtors}
Prior to proving Theorem~\ref{thm:sysfree4} of the Introduction, we
establish freedom modulo torsion for an arbitrary simplicial complex
of dimension~$4$, as follows.

\begin{theorem}
\label{thm:sys4}
Every compact, triangulable $4$-dimensional CW-complex $X$ is
$2$-systolically free (modulo torsion).
\end{theorem}

\begin{proof}
The map $X\to K(\Z^b,2)$ of Lemma \ref{cohotrick}, followed by a
homotopy equivalence $K(\Z^b,2)\simeq (\CP^{\infty})^{\times b}$ gives a map
$i:X\to (\CP^{\infty})^{\times b}$.  Up to homotopy, we may assume that the
image of $i$ lies in the $4$-skeleton $K_{b}=K(\Z^b,2)^{(4)}$.  By
construction, the map $i$ induces an isomorphism in $2$-dimensional
homology with integer coefficients.  By Lemma \ref{lem:longcyl2}, the
systolic freedom of $X$ reduces to that of the $4$-dimensional complex
$K_{b}$.

The map to $K_b$ constructed in the proof of the theorem is only an
isomorphism in rational $2$-dimensional homology, so we have no
control over areas of $2$-cycles defining torsion classes.

The CW-complex $K_{b}$ contains no $3$-cells.  Indeed, it is obtained
from the bouquet $\bigvee^b S^2_r$ of $b$ copies of $S^2$ by attaching
$4$-cells along all the Whitehead products $\frac{1}{2}[e_r,e_r]$ and
$[e_r,e_s]$, where $e_r$ is the fundamental class of $S^2_r$
(\cf paragraph~\ref{subsec:hm}).  Therefore, the freedom of $K_{b}$ 
reduces to that of the closures of its $4$-cells, by
Lemma~\ref{lem:closures}, since the isoperimetric inequality for
small $1$-cycles obviously holds for bouquets of spheres.

Each such closure is homeomorphic to either $\CP^2$ or 
$S^2 \times S^2$.  We are thus left with proving the 
systolic freedom of these two manifolds.

The freedom of $\CP^2$ reduces to that of the product of spheres as
follows.  Notice that there is a degree~$4$ map from $\CP^2$ to $P=S^2
\cup_{[e,e]} D^4$, where $e$ is the fundamental class of $S^2$.  Now
$P$ is homotopy equivalent to $W = S^2\times S^2 \cup_{a+b} D^3$,
where $a$ and $b$ are the fundamental classes of the two factors.
Thus $\CP^2$ admits a meromorphic map to $S^2\times S^2$, and we
invoke Proposition~\ref{prop:xyfree}.

The product of spheres is systolically free by
Corollary~\ref{cor:s2s2mero}, Lemma~\ref{lem:frees2t2}, and
Lemma~\ref{lem:longcyl2}.
\end{proof}

\begin{corollary}
\label{cor:fourmfd}
Every compact, smooth $4$-manifold $X$ admits metrics of arbitrarily
small volume, with the following property: every orientable, immersed
surface of smaller than unit area, defines a torsion class in
$H_2(X,\Z)$.
\end{corollary}

\begin{proof}
This is immediate from the theorem.
\end{proof}
Note that the manifold $X$ in the above corollary may be non-orientable.

\subsection{Orientable case}
\label{subsec:orientfree}
We now establish the systolic $2$-freedom of orientable $4$-manifolds.

\begin{proof}[Proof of Theorem~\ref{thm:sysfree4}]
Let $X$ be a compact, orientable, smooth $4$-mani\-fold.  By
Corollary~\ref{simplify} and Proposition~\ref{prop:xyfree}, the
freedom of $X$ reduces to that of a simply-connected manifold $Y$.
Since $Y$ is simply-connected, the group $H_2(Y;\Z)$ is free abelian.
Hence the $2$-freedom of $Y$ follows from Theorem~\ref{thm:sys4}.
\end{proof}

\subsection{The case of $2$-systolic freedom over $\Z_2$}
\label{subsec:mod2}
The notion of freedom with coefficients modulo $2$ was introduced in
Definition~\ref{def:z2free}.  We now reduce the question of
$2$-systolic freedom over $\Z_2$ of arbitrary $4$-manifolds to that of
just two manifolds: $S^2\times S^2$ and $\CP^2$.

\begin{proof}[Proof of Theorem \ref{nonorient}]
The orientable case is reduced to the case when $X$ is simply
connected, as in the proof of the previous theorem.  The map to
$K(\pi_2(X),2)$ induces an isomorphism in $2$-dimensional homology,
whether with integer or mod~$2$ coefficients.  Therefore the problem
is further reduced to the freedom over $\Z_2$ of this
Eilenberg-MacLane space.  The long cylinder construction of
Appendix~\ref{sec:freepull} works equally well with
$\Z_2$-coefficients.  Thus we may carve up the space $K(\pi_2(X),2)$
as above to reduce the problem to $\CP^2$ and $S^2\times S^2$.

The reduction of $\CP^2$ to $S^2 \times S^2$ as in the proof of
Theorem~\ref{thm:sys4} does not work here, as the map
$H_2(\CP^2;\Z_2)\to H_2(S^2\cup_{[e,e]}D^4;\Z_2)$ is not injective.
\end{proof}

\section{Systolic $2$-freedom in arbitrary dimension}
\label{sec:twofreedom}
The absence of odd cells in a suitable decomposition of the
Eilen\-berg-MacLane space $K(\Z,2)=\CP^{\infty}$ is the key to the
proof of Theorem~\ref{thm:sys4}.

The method employed in that proof can be used to improve
Theorem~\ref{thm:introbeyond} for $m=2$.  We illustrate this
by reducing the systolic $2$-freedom of all manifolds $X$ with
torsion-free $H_2(X,\Z)$ to that of a particularly simple
list of manifolds.

\begin{theorem}
\label{thm:red2free}
The following two statements are equivalent:
\begin{romenum}
\item All compact, smooth manifolds $X$ with torsion-free $H_2(X)$ are
$2$-systolically free;
\item For each $n\geq 2$, the manifold $\CP^n$ is $2$-systolically
free.
\end{romenum}
\end{theorem}

\begin{proof}
The $2$-freedom of each projective space would imply that of the
product of arbitrarily many factors, as in Lemma~\ref{further}.  The
$2$-freedom of a product of several copies of $\CP^1$ follows from
Corollary~\ref{cor:s2s2mero}.  But each cell closure in the standard
simplicial structure of $K(\Z^b,2)$ is such a product, proving that
any finite skeleton would also be $2$-free.  We map $X$ to $K(\Z^b,2)$
as in Lemma~\ref{cohotrick} and apply the pull-back
Lemma~\ref{lem:longcyl2}.
\end{proof}

\section{Odd-dimensional freedom}
\label{subsec:oddsys}
Note that the systolic $m$-freedom of $X$ when $m > \frac{n}{2}$
follows from the $m$-freedom of a bouquet of $m$-spheres viewed as an
$n$-dimensional CW-complex, by Corollary~\ref{pathological}.  Thus the
interesting case is $n>2m$.

\begin{theorem}
\label{thm:oddsys}
Let $n$ and $m$ be integers satisfying $3\leq m <n$, where $m$ is odd.
Then every $n$-dimensional manifold $X$ is $m$-systolically free
(modulo torsion).
\end{theorem}

\begin{proof}
The starting point is again the map from the $(m+1)$-skeleton of $X$
to the bouquet of spheres, as in Lemma~\ref{cohotrick}.  For odd $m$,
the only non-trivial Whitehead products are the `mixed' ones.  Thus,
high-degree self-maps of the sphere allow us to map $X$ to the product
$(S^m)^{\times b}$ of $b=b_m(X)$ copies of $S^m$, as in the proof of
Lemma~\ref{lem:prebouquet}.

Here we rely upon the existence of self-maps inducing the zero
homomorphism in every given homotopy group of the sphere.  Indeed,
according to D.~Sullivan~\cite{Sv}, a self-map of $S^m$ of degree $d$
induces a nilpotent map in the $d$-torsion of $\pi_j(S^m)$, for every
$j>0$.

Now if $m$ does not divide $n$, the $n$-skeleton of $(S^m)^{\times b}$
coincides with the $(n-1)$-skeleton, in which case we actually get
(singular) metrics on $X$ of vanishing $n$-volume.

If $m$ divides $n$, the $n$-skeleton of $(S^m)^{\times b}$ contains no
$(n-1)$-cells, and we use the carving-up procedure of
\S\ref{sec:carveup} to reduce to products of spheres.
\end{proof}

\appendix
\section{Systoles of complementary dimensions}
\label{sec:compdim}
To describe phenomena along the lines of Gromov's example, it is
convenient to introduce the following terminology.

\begin{definition}
\label{def:sysfree}
A compact, smooth, $(p+q)$-dimensional manifold $X$ is called
{\em systolically $(p,q)$-free} if
\[
\inf_{g}\frac{\vol_{p+q}(g)}{\sys_p(g) \sys_q(g)}= 0,
\]
where the infimum is taken over all metrics $g$ on $X$.
\end{definition}

We present a proof of $(1,n-1)$-freedom, originally obtained by the
first author in collaboration with I.~Babenko in \cite{BaK}.

\begin{theorem}
\label{thm:compfree}
Every compact, orientable, smooth manifold of dimension $n\geq 3$ is
$(1,n-1)$-free.
\end{theorem}

\begin{proof}
We obtain systolically free families of metrics on the manifold $X^n$
by a direct geometric construction.  The idea is to introduce a local
version of Gromov's example (\cf \S\ref{subsec:gromov}), \ie, a metric
on a manifold with boundary which can be glued into any manifold to
ensure systolic freedom.  Here the matrix of Gromov's example:
\[
\begin{bmatrix}
1&-j\\ -j & 1+j^2
\end{bmatrix}
\begin{matrix}
\leftarrow  z \\
\leftarrow  b
\end{matrix}
\]
for a given $j$, is replaced by the matrix
\[
\begin{bmatrix}
1&-x\\ -x & 1+x^2
\end{bmatrix}
\begin{matrix}
\leftarrow  z \\
\leftarrow  y
\end{matrix}
\]
where the coefficient $x$ varies between $0$ and $j$.

Let $C$ be a union of closed curves which form a basis for $H_1(X; \Q)$.
Its tubular neighborhood is diffeomorphic to $C\times D^{n-1}$. Let
$K\subset D^{n-1}$ be a codimension~$2$ submanifold with trivial
normal bundle (\eg, the $(n-3)$-sphere).  The boundary of a tubular
neighborhood of $C\times K$ is diffeomorphic to the hypersurface
$\Sigma:= C\times K\times T^1$, where $T^1$ is a circle.  Our
construction is local in a neighborhood of $\Sigma$.  Choose a fixed
metric on $X$ satisfying the following four properties:

\begin{enumerate}

\item\label{1}
It is a direct product in a neighborhood of $\Sigma$.

\item\label{2}
The hypersurface $\Sigma$ is a metric direct product of
the three factors $C$, $K$, and $T^1$.

\item\label{3}
Each connected component of $C$ is a circle of length $1$.

\item\label{4}
The circle $T^1$ has length $1$.

\end{enumerate}

We now cut $X$ open along $\Sigma$ and insert suitable `cylinders'
$\Sigma\times I$, indexed by $j\in \N$, resulting in a sequence of
metrics $g_j$ on $X$.  These `cylinders' are {\em not} metric
products.  They have the following properties.
\begin{enumerate}
\setcounter{enumi}{4}
\item\label{5}
The projection $\Sigma\times I\rightarrow I$ is a Riemannian
submersion over an interval $I=[0,2j]$, where the interesting behavior
is exhibited when $j$ grows without bound.

\item\label{6}
The metric at the endpoints $0$ and $2j$ agrees with that of $\Sigma$.

\item\label{7} The $(n-3)$-dimensional manifold $K$ has a fixed metric
independent of $j$, and is a direct summand in a metric product.

\item\label{8}
For integer values of $x\in I$, each connected component of
$C\times T^1\times \{x\}$ is isometric to a standard unit
square torus.

\item\label{9}
The metric on $C\times T^1\times [0,2j]$ is the `double' of the
metric on $C\times T^1\times [0,j]$, in the sense that $C\times
T^1\times \{x\}$ and $C\times T^1\times \{2j-x\}$ have identical
metrics.

\end{enumerate}

Now let $N$ be a non-trivial $(n-1)$-cycle of $X$.  By Poincar\'{e}
duality, the cycle $N$ has non-zero intersection number with one of
the connected components, $C_i$, of $C$.  Hence $N$ induces a
non-trivial relative cycle in a neighborhood of $C_i$.  From now on we
will denote this component by $C$.

Note that the volume of $(X,g_j)$ grows linearly in $j$.  We will
obtain a lower bound for the $(n-1)$-volume of $N$, and therefore for
the $(n-1)$-systole of $g_j$, which grows faster than the volume of
$g_j$.  Meanwhile, the $1$-systole is bounded from below uniformly in
$j$.  The theorem now follows from the properties of suitable metrics
on $T^2\times I$ constructed below.

Our technique is calibration by the $(n-1)$-form $\alpha\wedge
\mu_K=*dz\wedge\mu_K$, where $\mu_K$ is the volume form of $K$.  Here
the $2$-form $\alpha= *dz$ provides the lower bound for the area of
the relative $2$-cycle $M$ in Lemma~\ref{lem:hodge} below.
\end{proof}

The metric on $C\times T^1\times [0,2j]$ is {\em not} a direct sum.
Consider the subinterval $[0,j]\subset[0,2j]$.  Then $C\times T^1\times
[0,j]$ can be thought of as a fundamental domain for a $j$-fold cover of
a non-trivial torus bundle over the circle, defining either the
$\operatorname{NIL}$ or $\operatorname{SOL}$ geometries (used
in \cite{BeK} and \cite{P}, respectively).
We will present a description valid for both approaches.

What makes these metrics systolically interesting are the following
properties.  There is an orientable surface
$M^2\subset C\times T^1\times [0,2j]$ and
a $2$-form $\alpha$ on $C\times T^1\times [0,2j]$ such that
\begin{alphenum}

\item\label{ia} $M$ and $C$ have intersection number equal to $1$.

\item\label{ib} $\length(M\cap T^2_x)\geq x$ for all $x\in [0,j]$,
where $T^2_x= C\times T^1\times \{x\}$.

\item\label{ic}
$\alpha(e_1,e_2) = |\alpha(p)|$ for every $p\in M$, where
$(e_1,e_2)$ is an orthonormal basis of $T_p M$ (\ie, the $2$-form is
maximal in the direction tangent to $M$).

\item\label{id}
The $2$-form $\phi(x) \alpha$ is closed, for any function $\phi$.
\end{alphenum}

More specifically, we have the following.

\begin{lemma}
\label{lem:hodge}
Equip the manifold $Y= T^2\times I$ with metrics $g_j$ defined by
\[
g_j(x,y,z)= h(\hat x) (y,z) + dx^2,
\]
where $x\in I=[0,2j]$, $\hat x= \min (x,2j-x)$, $T^2$ is the quotient
of the $(y,z)$-plane by the integer lattice, and $h(x) (y,z) =
(dz-xdy)^2 + dy^2$ defines a metric on the $2$-torus $T^2\times \{x\}$
(see Figure~\ref{fig:thinpars}).
Then $Y$ has properties \textup{(\ref{ia})--(\ref{id})}.
\end{lemma}

\setlength{\unitlength}{0.95cm}
\begin{figure}
\begin{picture}(11.5,4.2)(0,0.5)
\put(0,1){\line(1,0){11.5}}
\multiput(0,1)(1.5,0){3}{\line(0,1){0.2}}
\put(5.5,1){\line(0,1){0.2}}
\multiput(8.5,1)(1.5,0){3}{\line(0,1){0.2}}
\put(0,2){\framebox(0.5,0.5)}
\multiput(1.5,2)(0.5,0.5){2}{\line(0,1){0.5}}
\multiput(1.5,2)(0,0.5){2}{\line(1,1){0.5}}
\multiput(3,2)(0.5,1){2}{\line(0,1){0.5}}
\multiput(3,2)(0,0.5){2}{\line(1,2){0.5}}
\multiput(5.5,2)(0.5,2){2}{\line(0,1){0.5}}
\multiput(5.5,2)(0,0.5){2}{\line(1,4){0.5}}
\multiput(8.5,2)(-0.5,1){2}{\line(0,1){0.5}}
\multiput(8,3)(0,0.5){2}{\line(1,-2){0.5}}
\multiput(9.5,2.5)(0.5,-0.5){2}{\line(0,1){0.5}}
\multiput(9.5,2.5)(0,0.5){2}{\line(1,-1){0.5}}
\put(11,2){\framebox(0.5,0.5)}
\put(0,0.7){\makebox(0,0){\footnotesize{$0$}}}
\put(1.5,0.7){\makebox(0,0){\footnotesize{$1$}}}
\put(3,0.7){\makebox(0,0){\footnotesize{$2$}}}
\put(5.5,0.7){\makebox(0,0){\footnotesize{$j$}}}
\put(11.5,0.7){\makebox(0,0){\footnotesize{$2j$}}}
\put(-0.15,2.25){\makebox(0,0){\scriptsize{$1$}}}
\put(0.25,1.8){\makebox(0,0){\scriptsize{$1$}}}
\put(1.35,2.25){\makebox(0,0){\scriptsize{$1$}}}
\put(1.88,2.1){\makebox(0,0){\scriptsize{$\sqrt{2}$}}}
\put(2.85,2.25){\makebox(0,0){\scriptsize{$1$}}}
\put(3.48,2.45){\makebox(0,0){\scriptsize{$\sqrt{5}$}}}
\put(5.35,2.25){\makebox(0,0){\scriptsize{$1$}}}
\put(6.43,3.2){\makebox(0,0){\scriptsize{$\sqrt{1+j^2}$}}}
\end{picture}
\caption{\textsf{Family of parallelograms defining metrics $g_j$ on $T^2
\times I$}}
\label{fig:thinpars}
\end{figure}
\setlength{\unitlength}{1cm}

\begin{proof}
Define the surface $M$ to be the cylinder $T^1\times I\subset
\R^2/\Z$, an open subset of the quotient of the $(y,x)$-plane by unit
translation in the $y$-coordinate.  The curve $C$ is the projection of
the $z$-axis.  Set $\alpha=* dz$, the Hodge star of the coordinate
$1$-form $dz$, for the metric $g_j$.

With respect to the metric $g_j$, the forms $dx$, $dy$, and $dz-xdy$
form an orthonormal basis of $T_{(x,y,z)}^* Y$.
We will need the following calculation of the pointwise norm of the
$1$-form $dz$:
\[
\abs{dz}^2=\abs{dz-xdy}^2 +\abs{xdy}^2=1+x^2.
\]

Now given any surface $M'\subset(Y,g_j)$ in the relative homology
class of $M$ in $H_2(Y,\partial Y)$, we have by Stokes' theorem and
properties (\ref{ic}) and (\ref{id}),
\[
\area(M')\geq\int_{M'}\phi(x)*\!dz=
\int_{M}\phi(x)*\!dz\geq\int_1^{j-1}\sqrt{1+x^2}dx\sim{j^2},
\]
for any function $\phi(x)$ with support in $[0,j]$, which is dominated
by the function $\abs{dz}^{-1}$ and coincides with it on $[1,j-1]$.  This
lower bound for the relative $2$-systole is the source of all freedom.
\end{proof}

\begin{remark}
\label{mod2}
The metric $g_j$ does admit a $\bmod 2$ relative $2$-cycle in the same
class as $M$, whose area grows linearly in $j$.  Namely, consider the
$2$-chain $a=T^1\times [0,2]\subset C\times T^1\times [0,2j]$.  Let
the $2$-chain $b\subset T^2\times \{2\}$ be defined by the projection
to the torus of the triangle of base 1 (the circle $T^1$) and altitude
$2$ (twice the curve $C$).  The boundary of $b$ modulo $2$ consists of
two circles: image of the base and image of the hypothenuse.  Let
$c=a+b$.  Let $h:T^2\times \{0\}\to T^2\times \{1\}$ be the glueing
homeomorphism of the $\operatorname{NIL}$ bundle which provided the
starting point for our construction.  Here $h$ is represented by the
matrix $\bigl(\begin{smallmatrix}1&1\\0&1\end{smallmatrix}\bigr)$ and
is an isometry Let with respect to the metric $g_j$.
\[
d=c+ h^2(c)+ h^4(c) + h^6(c) +\cdots +h^{j-2}(c)
\]
(here we are assuming that $j$ is even).  Then $d$ is a $\bmod\: 2$
relative $2$-cycle of $T^2\times [0,j]$.  Doubling it as in
item~\ref{9} above, we obtain the desired relative cycle, representing
the generator of $H_2(T^2\times I, \partial(T^2\times I);\Z_2)$.
\end{remark}

For the inequalities involving systoles of complementary dimensions
$k$ and $n-k$, the existing results on freedom depend on the
divisibility of $k$ by $4$.

\begin{theorem}[\cite{BaK}]
\label{thm:divby4}
Let $X^n$ be a compact, orientable, smooth $n$-dimensional manifold.
Assume $X$ is $(k-1)$-connected, where $k<\frac{n}{2}$ and $k$ is not
divisible by $4$.  Then $X$ is $(k,n-k)$-systolically free.
\end{theorem}

\section{Freedom pulled back}
\label{sec:freepull}
Our main goal here is to prove Lemma~\ref{lem:freep}, which follows
from the three lemmas below.

\begin{lemma}[\cite{BaK}]
\label{lem:longcyl1}
Let $Y$ be a triangulable CW-complex of dimension~$n$.  Let
$W=Y\cup_{f} e^{k}$ be a complex obtained from $Y$ by attaching a
single cell of dimension $k\leq n-1$, where the attaching map $f$ is
triangulable.  Let $1\leq q\leq n-1$.  Then the $q$-systolic freedom
of $Y$ implies that of $W$.  The same statement is true for
$(q,n-q)$-freedom, etc.
\end{lemma}

\begin{proof}
Here the volume of an $n$-dimensional triangulable complex is by
definition the sum of volumes of all cells of maximal dimension.

The idea is to insert a long $k$-dimensional cylinder in such a way
that a $q$-cycle of volume comparable to the $q$-systole of $Y$ would
necessarily have a `narrow place' somewhere along the cylinder, by
virtue of the coarea inequality.  We cut the cycle into two pieces at
the narrow place, and fill the cut with a small $q$-chain (using the
isoperimetric inequality for small $(q-1)$-cycles) to make both pieces
into cycles.  The lemma now follows from the fact that the cylinder,
being of positive codimension, does not contribute to top-dimensional
volume.

By the `long cylinder' metric we mean the following.  Let $I=[0,\ell]$,
with $\ell\gg 1$ to be determined.  Let $S=\Image(f)$ be the image
of the attaching map.  We may assume that $S$ is a subcomplex of $Y$.

Let $\Cyl_f=Y\cup_{f\times\{0\}}(S^{k-1}\times I)$ be the mapping
cylinder of $f$ and $C_f=\Cyl_f\cup_{\id\times\{\ell\}}D^k$ the
mapping cone, where $D^k$ is a cell of dimension $k\leq n-1$.

Let $g$ be a metric on $Y$.  By scaling, we may assume that
$\sys_q(g)\geq 1$.

Let $h_0$ be the restriction of $g$ to $S\subset Y$.
Let $(S^{k-1},h_1)$ be a round sphere of sufficiently
large radius $r\geq 1$ so that the triangulable map
$f:(S^{k-1},h_1)\to (S,h_0)$ is distance-decreasing,
or roughly, $h_1\geq h_0$.  Let $D^k$ be a round hemisphere
of the same radius as $(S^{k-1},h_1)$.

We endow the cylinder $S^{k-1}\times I$ with the metric $(1-x)h_0+x
h_1+dx^2$ for $0\leq x\leq 1$ and $h_1+dx^2$ for $1\leq x\leq \ell.$
We thus obtain a metric on the complex $W=C_f$.  Denote the resulting
metric space by $W(g,\ell)$.

Note that mapping cylinder $\Cyl_f\subset W(g,\ell)$ admits a
distance-decreasing projection to $(Y,g_j)$ by construction.

Let $p:W(g,\ell)\rightarrow I$ be the map extending to $W$ the
projection to the second factor $S^{k-1}\times I\rightarrow I$ on the
cylinder, while $p(Y)=0$ and $p(D^k)=\ell$.  Let $z$ be a $q$-cycle in
$W$. The complex $W$ is not a manifold, but its subspace
$S^{k-1}\times I$ is a manifold and we can apply the coarea inequality
just in this part of $W$.  We obtain the inequality $\vol_q(z)\geq
\int_1^{\ell} \vol_{q-1}(z\cap p^{-1}(x))dx$.  Hence we can find a
regular value $x_0\in I$ such that the $(q-1)$-cycle $c= z\cap
p^{-1}(x_0)$ satisfies
\begin{equation}
\label{eq:volz}
\vol_{q-1}(c)= \vol_{q-1}(z\cap p^{-1}(x_0))\leq
\tfrac{1}{\ell-1}\vol_q(z).
\end{equation}
Note that in the case $q=k-1$, the cycle $c$ is empty.

Let us now show that if $Y$ admits a systolically free sequence of
metrics $g_j$, then so does $W$.  (Note that we are not assuming that
the inclusion of $Y$ in $C_f$ induces a surjective map in
$q$-dimensional homology.)

Choose $\ell=\ell(j)\geq \sys_q(g_j)\geq 1$.  We would like to show
that the $q$-volume of non-trivial $q$-cycles in $W(g_j,\ell)$ is
comparable to $\sys_q(g_j)$.  Suppose, on the contrary, that a
sequence of cycles $z_j$ in $W(g_j,\ell)$ satisfies
\[
\vol_q(z_j)=o(\sys_q(g_j)).
\]
Letting $c_j=z_j\cap p^{-1}(x_0)$, for suitable $x_0$ depending on $j$
as in \eqref{eq:volz}, we obtain $\lim_{j\rightarrow\infty}
{\vol_{q-1}}(c_j)=0$.

We now appeal to the isoperimetric inequality for small cycles (\cf
\cite{G1}, Sublemma $3.4.B'$), applied to the cross section
$S^{k-1}\times \{x_0\}$ of the cylinder, to obtain a $q$-chain
$B^q_j\subset S^{k-1}$ satisfying
$$\vol_q(B^q_j)\leq K \vol_{q-1}(c_j)^\frac{q}{q-1}.$$ In other words,
the $(q-1)$-cycle $c_j$ can be filled by a chain whose volume also
tends to $0$.  Clearly, the constant $K$ in the inequality can be
chosen independent of the radius $r\geq 1$ of $S^{k-1}$.  Let
\[
a_j= (z_j\cap d^{-1}([0,x_0])) - B^q_j
\]
and $b_j=z_j-a_j$.  Note that $[b_j]=0$ in $W$ and so
$[a_j]=[z_j]\not=0$.  The cycle $a_j$ lies in the mapping cylinder
which admits a distance-decreasing projection to $(Y,g_j)$.

Hence $\vol_{q}(a_j)\geq \sys_{q}(g_j)$ and so $\vol_{q}(z_j)\geq
\vol_{q}(a_j)-\vol_{q}(B^q_j)\geq \sys_q(g_j)-o(1)$.  This shows that
the systoles of $W$ are not significantly diminished compared to those
of $Y$.
\end{proof}

\begin{lemma}
\label{lem:longcyl2}
Let $X$ and $W$ be triangulable CW-complexes of dimension~$n$. Suppose
there is a map $\phi: X \to W$ inducing a monomorphism on $H_q$.  Then
the $q$-systolic freedom of $W$ implies that of $X$.  The same
statement holds for $(q,n-q)$-freedom, etc.
\end{lemma}

\begin{proof}
Choose a simplicial structure on $W$.  By the cellular approximation
theorem, a continuous map from $X$ to $W$ is homotopic to a simplicial
map.  As in \cite{Ba}, we can replace it by a map which has the
following property with respect to suitable triangulations of $X$ and
$W$: on each simplex of $X$, it is either a diffeomorphism onto its
image or the collapse onto a wall of positive codimension.  Let $p$ be
the maximal number of $n$-simplices of $X$ mapping diffeomorphically
to an $n$-simplex of $W$.  Since cells of dimension $\leq (n-1)$ do
not contribute to $n$-dimensional volume, the pullback of the metric
on $W$ is a positive quadratic form on $X$ whose $n$-volume is at most
$p$ times that of $W$.  This form is piecewise smooth and satisfies
natural compatibility conditions along the common face of each pair of
simplices.
\end{proof}

\begin{lemma}
\label{lem:longcyl3}
If a smooth compact $n$-manifold $X$ admits systolically free
piecewise smooth metrics, then it also admits systolically free smooth
metrics.
\end{lemma}

\begin{proof}
To construct a smooth metric from a piecewise smooth one, we proceed
as in \cite{Ba}.  Given a piecewise smooth metric $g$, compatible
along the common face of each pair of adjacent simplices, we choose a
smooth metric $h$ on $X$ such that $h > g$ at every point (in the
sense of lengths of all tangent vectors).  Let $N$ be a regular
neighborhood of small volume of the $(n-1)$-skeleton of the
triangulation.  Choose an open cover of $X$ consisting of $N$ and the
interiors $U_i$ of all $n$-simplices.  Using a partition of unity
subordinate to this cover, we patch together the metrics $g|_{U_i}$
and $h|_N$.  The new metric dominates $g$ for each tangent vector to
$M$.  In particular, the volume of a cycle is not decreased.
Meanwhile, $n$-dimensional volume is increased no more than the volume
of the regular neighborhood.

The piecewise smooth metric on $X$ may {\it a priori} not be
compatible with its smooth structure, since the triangulation may not
be smooth.  To clarify this point, denote the triangulation by $s$,
and the piecewise smooth metric by $g$.  Consider a smooth
triangulation $s'$, and approximate the identity map of $X$ by a map
$\phi$, simplicial with respect to $s'$ and $s$.  The pullback metric
$g'=\phi^*(g)$ is then adapted to the smooth triangulation $s'$, so we
may apply to it the argument with the regular neighborhood $N$.

We have thus obtained a smooth positive form on $X$.  We make it
definite without significantly increasing its volume by adding a small
multiple of a positive definite form.  The lower bounds for the
$q$-systole are immediate from the injectivity of the map
$X\rightarrow W$ on $H_q$.
\end{proof}

We now establish the systolic freedom of a model space by a `long
cylinder' argument.

\begin{lemma}
\label{lem:longcyl4}
Let $X$ and $Y$ be $q$-systolically free $n$-complexes, and let $S$ be
an $(n-2)$-complex which is simultaneously a subcomplex of $X^{(n-2)}$
and $Y^{(n-2)}$.  Let $I$ be an interval.  Then the complex
\[
W=X\cup (S\times I)\cup Y
\]
is also $q$-systolically free.
\end{lemma}

\begin{proof}
Let $\epsilon>0$, and let $(X,g)$ and $(Y,h)$ be metrics satisfying
\[
\frac{\vol(g)^{\frac{q}{n}}}{\sys_q (g)}<\epsilon\quad \text{ and }\quad
\frac{\vol(h)^{\frac{q}{n}}}{\sys_q (h)}<\epsilon.
\]
We will construct a metric on $M$ satisfying a similar inequality for
$2\epsilon$. For this purpose, assume that $g$ and $h$ are scaled to
the same value of the $q$-systole.

Let $\mu$ be a metric on $S$ dominating both (the pullbacks of) $g$
and $h$.  Let $I=[0,\ell]$.  We define a metric on $W$ to be equal to
$\mu\oplus dt^2$, when $t\in[1,\ell-1]\subset I$.  Near the
extremities of $I$, we patch the metric to make a continuous
transition to $g$ and $h$, as in the proof of
Lemma~\ref{lem:longcyl1}.  The systolic freedom of $W$ follows by
applying the coarea inequality as before.  Namely, a non-trivial
$q$-cycle in $W$ can be cut at a narrow place and decomposed, by
filling the cut, into two pieces, one of which admits a
distance-decreasing projection to $(X,g)$, and the other to $(Y,h)$.
One of the pieces must be non-trivial (since they add up to a
non-trivial cycle), hence bounded below by $\sys_q(X)=\sys_q(Y)$.  The
contribution from the filling cycle is negligible by the isoperimetric
inequality for small cycles.
\end{proof}

\section{A Mayer-Vietoris argument}
\label{sec:mv}
\begin{lemma}
\label{lem:circsurg}
Let $X$ be a closed, orientable, smooth $4$-manifold.  Let $C\subset X$ 
be a smoothly embedded closed curve in $X$, representing an element
in the torsion-free part of $H_1(X)$.  Let $Y$ be the result of
surgery along $C$.  Let $Z$ be the cobordism determined by the
surgery, and let $i:X\to Z$ and $j:Y\to Z$ be the inclusions.  There
is then an isomorphism $\psi:H_2(X)\to H_2(Y)$ such that 
$j_*\circ \psi = i_*$.
\end{lemma}

\begin{proof}
Let $k:S^1\to X$ be an embedding with $k(S^1)=C$.
By assumption, $k_*: H_1(S^1)\to H_1(X)$ is a monomorphism.
Let $X_{-}$ be a neighborhood of $S^1$ (homeomorphic to
$S^1\times D^3$), and $X_{+}$ the closure of its complement.
Let $i^{\pm}: X_{\pm}\to X$
and $k^{\pm}:S^1\times S^2\to X_{\pm}$ be the inclusions.
Consider the following fragment of the Mayer-Vietoris sequence:
\[
\begin{split}
H_3(X) \xrightarrow{\partial} H_2(S^1\times S^2)
\rightarrow H_2(X_{-})\oplus H_2(X_{+})\xrightarrow{(i^{-}_{*},i^{+}_{*})}
H_2(X) \rightarrow  \\
H_1(S^1\times S^2) \xrightarrow{k^{-}_{*}+k^{+}_{*}}
H_1(X_{-})\oplus H_1(X_{+})\xrightarrow{(i^{-}_{*},i^{+}_{*})}
H_1(X)\to 0.
\end{split}
\]
We have: $H_2(X_{-})=0$, the map $\partial$ is surjective
(under Poincar\'{e} duality and universal coefficients, it is
the dual of $k_*: H_1(S^1)\to H_1(X)$, which is injective),
and $k^{-}_{*}+k^{+}_{*}$ is injective, with image $H_1(X_{-})$.
Thus, the maps $i^{+}_{*}: H_1(X_{+})\to H_1(X)$ and
$i^{+}_{*}: H_2(X_{+})\to H_2(X)$ are isomorphisms.

The surgery replacing $X_{-}$ by $D^2\times S^2$ yields
$Y=D^2\times S^2 \cup X_{+}$.  Let $j^{-}:D^2\times S^2\to Y$,
$j^{+}:X_{+}\to Y$, and $\ell:S^1\times S^2\to D^2\times S^2$
be the inclusions.  The corresponding Mayer-Vietoris sequence is
\[
\begin{split}
H_2(S^1\times S^2) \xrightarrow{\ell_* + k^{+}_*}
H_2(D^2\times S^2)\oplus H_2(X_+) \xrightarrow{(j^{-}_{*},j^{+}_{*})}
H_2(Y) \rightarrow \\
H_1(S^1\times S^2)\xrightarrow{\ell_* + k^{+}_*}
H_1(D^2\times S^2)\oplus H_1(X_{+}).
\end{split}
\]
We have: $H_1(D^2\times S^2)=0$, the map
$\ell_*: H_2(S^1\times S^2) \to H_2(D^2\times S^2)$
is an isomorphism, and $k^{+}_*: H_1(S^1\times S^2)\to H_1(X_{+})$
is injective (since $k_*$ is injective, $i^{+}_*$ is an isomorphism, and
$i^{+}_*\circ k^{+}_*=k_*$). Thus, $j^{+}_*: H_2(X_{+})\to H_2(Y)$
is an isomorphism.

The desired isomorphism is obtained by combining the
two isomorphisms above:
\[
\psi=j^{+}_{*}\circ (i^{+}_{*})^{-1}:H_2(Y)\to H_2(X).
\]
It remains to verify that $j_*\circ \psi = i_*$.  Recall that, by
definition, $Z= X\times I \cup h^2$, where the handle $h^2=D^2\times
D^3$ is attached along $X_{-}=S^1\times D^3$.  We thus have a
diffeomorphism $Z\cong X_{+}\times I \cup (X_{-}\times I\cup h^2)$.
Under this decomposition, the inclusions $i\circ i^{+}$ and $j\circ
j^{+}$ of $X_{+}$ into the boundary components of $Z$ correspond to
the inclusions $\iota_0$ and $\iota_1$ of $X_{+}$ into the boundary
components of $X_{+}\times I$.  The identity map $X_{+}\times I \to
X_{+}\times I$ is a homotopy $\iota_0\simeq \iota_1$.  This implies
$i\circ i^{+}\simeq j\circ j^{+}$, hence $i_*\circ i^{+}_*=j_*\circ
j^{+}_*$, and the conclusion follows.
\end{proof}

\section{Spinning}
\label{sec:spin}
The process of spinning was introduced in knot theory by E.~Artin, and
was extended to closed manifolds by S.~Cappell~\cite{C}.

Let $M^n$ be a compact, simply-connected, smooth manifold of dimension
$n\ge 3$.
The result of performing surgery on $S^1\times M$ along $S^1\times *$ is
called the {\em spin} of $M$.  There is a choice of framing involved in the
surgery, corresponding to $\pi_1(\SO(n))=\Z_2$.  The spin obtained by
surgery with trivial framing is denoted by $\sigma(M)$; the `twisted'
spin, by $\sigma'(M)$.

Note that $\sigma(M)=\partial (D^2 \times M_{\circ})$, where
$M_{\circ}=M^n\setminus \interior D^n$ is a punctured copy of $M$.
Indeed, $\sigma(M)=S^1\times M^n\setminus (S^1\times \interior D^n)
\cup D^2\times S^{n-1} = S^1\times M_{\circ}\cup D^2\times S^{n-1}
=\partial D^2 \times M_{\circ} \cup D^2\times \partial M_{\circ}
=\partial (D^2 \times M_{\circ})$.

If $M$ admits an effective circle action with non-empty fixed point
set, then it is easily shown that $\sigma'(M)\cong \sigma(M)$.
In general, though, the two spins of $M$ are not even
homotopy equivalent.

\begin{lemma}[\cite{Su}]
\label{lem:sigma}
Both spins of $S^1\times S^{n-1}$ can be identified with the
connected sum $S^1 \times S^n \# S^{n-1} \times S^2$.
\end{lemma}

\begin{proof}  We have $\sigma(S^1 \times S^{n-1})=\partial W$,
where $W=D^2 \times (S^1\times S^{n-1})_{\circ}$.  The manifold
$W$ admits a handle  decomposition $W=D^{n+2}\cup h^1 \cup h^{n-1}$, where
the attaching spheres form a trivial link $S^0\cup S^{n-2} \subset S^{n+1}$.
Thus, $W$ can be identified with the boundary connected sum
$(D^{n+2}\cup h^1)\sharp (D^{n+2}\cup h^{n-1}) \cong
S^1\times D^{n+1}\sharp S^{n-1}\times D^3$.  Passing to boundaries,
we obtain
\begin{equation}
\label{eq:spindiff}
\sigma(S^1\times S^{n-1})\cong
S^1 \times S^n \# S^{n-1} \times S^2.
\end{equation}

Rotation in the second factor of $S^1\times S^{n-1}$ provides a
circle action with non-empty fixed point set.  Thus,
$\sigma'(S^1\times S^{n-1})\cong \sigma(S^1\times S^{n-1})$.
\end{proof}

\begin{corollary}
\label{cor:spin}
The manifold $S^2 \times S^{n-1}$ is obtained from
$S^1 \times S^1\times S^{n-1}$ by performing surgery
along the curves $\gamma_1=S^1\times *\times *$ and
$\gamma_2=*\times S^1\times *$.
\end{corollary}

\begin{proof}
First consider the case where both surgeries are done with trivial
framing.  Surgery along $\gamma_1$ yields the spin $\sigma(S^1 \times
S^{n-1})$.  Under the diffeomorphism in \eqref{eq:spindiff}, the image
of the curve $\gamma_2$ on the left side corresponds to the curve
$S^1\times *$ in the factor $S^1 \times S^{n}$ on the right side.
Surgery on $S^1 \times S^{n}$ along $S^1\times *$ yields
$\sigma(S^{n})=S^{n+1}$, and the corollary follows.

The other cases are similar.  It suffices to note that $S^n$ admits an
effective circle action with fixed points, and so $\sigma'(S^n)\cong
\sigma(S^n)$.
\end{proof}

\section{Whitehead products and Hilton's theorem}
\label{sec:wh}

The product of two spheres decomposes as $S^i_a \times S^j_b = (S^i_a
\vee S^j_b)\cup_{f} e^{i+j}$, where $f:S^{i+j-1}\to S^i\vee S^j$ is
the attaching map of the top cell, so that $[f]=0$ in $\pi_{i+j-1}
(S^i\times S^j)$.  This decomposition defines a $\Z$-bilinear pairing
\[
[\ ,\ ]:\pi_i(X)\times\pi_j(X) \to \pi_{i+j-1}(X)
\]
on the homotopy groups of an arbitrary space $X$.  Here the element
$[a,b]$ in the case of our bouquet is defined by the map $f$.  More
precisely, if $\a\in \pi_i(X)$ and $\b\in \pi_j(X)$, then $[\a,\b]$ is
the obstruction to extending the map $a\vee b: S^i\vee S^j \to X$ to
the whole of $S^i\times S^j$, where $a$ and $b$ are arbitrary
representatives of $\a$ and $\b$.

The above operation, which generalizes the usual commutator map in
$\pi_1$, is called the {\em Whitehead product}.  It is natural with
respect to continuous maps, and satisfies identities analogous to the
Lie bracket.  See \cite{Wh} for more details.

\begin{theorem}[Hilton~\cite{Hil}]
\label{subsec:hm}
The first unstable homotopy group of a bouquet of spheres is given by:
\[
\pi_{2m-1}(\vee S^m_r)=\oplus_r \pi_{2m-1} (S^m_r) \oplus
(\oplus_{r<s} \Z[e_r,e_s]).
\]
Meanwhile, $\pi_j(\vee S^m_r)=\oplus_r \pi_j (S^m_r)$ for
$m\leq j\leq 2m-2$.
\end{theorem}

For example, $\pi_3(S^2 \vee S^2)=\Z \oplus \Z \oplus \Z$.  The first
$\Z$-summand corresponds to $\pi_3(S^2)$, and is generated by the Hopf
map $h:S^3\to S^2$.  This generator is written unambiguously (by
virtue of the absence of torsion in this dimension) as
$h=\frac{1}{2}[a,a]$, where $a$ is the fundamental class of $S^2$.
Similarly, the second $\Z$-summand is generated by $\frac{1}{2}[b,b]$.
The last $\Z$-summand is generated by $[a,b]$.

\bibliographystyle{amsalpha}

\begin{thebibliography}{99}

\bibitem{Ba} I.~Babenko,
{\em Asymptotic invariants of smooth manifolds},
Russian Acad. Sci. Izv. Math. \textbf{41} (1993), 1--38.

\bibitem{BaK} I.~Babenko and M.~Katz, {\em Systolic freedom of
orientable manifolds}, Ann. Sci. \'{E}cole Norm. Sup. \textbf{31}
(1998); available at {\ttfamily math.DG/9707102.}

\bibitem{BKS} I.~Babenko, M.~Katz and A.~Suciu,
{\em Volumes, middle-dimensional systoles, and Whitehead products},
Math. Res. Lett. \textbf{5} (1998), 461--471.

\bibitem{BeK} L.~B\'{e}rard Bergery and M.~Katz,
{\em Intersystolic inequalities in dimension $3$},
Geom. Funct. Anal. \textbf{4} (1994), 621--632.

\bibitem{Ber1} M.~Berger,
{\em Du c\^{o}t\'{e} de chez Pu},
Ann. Sci. \'{E}cole Norm. Sup. \textbf{5} (1972), 1--44.

\bibitem{Ber2} \bysame,
{\em Systoles et applications selon Gromov},
S\'eminaire N.~Bourbaki, expos\'{e} 771, Ast\'{e}risque
\textbf{216} (1993), 279--310.

\bibitem{C}  S.~Cappell,
{\em Superspinning and knot complements},
In: Topology of Manifolds
(Proc. Inst. Univ. of Georgia, Athens, Ga., 1969),
Markham, Chicago, IL, 1970,  pp. 358--383.

\bibitem{E} B.~Eckmann,
{\em \"{U}ber die Homotopiegruppen von Gruppenr\"{a}umen},
Comment. Math. Helv. \textbf{14} (1942), 234--256.

\bibitem{Fe} H.~Federer,
{\em Real flat chains, cochains and variational problems},
Indiana Univ. Math. J. \textbf{24} (1974), 351--407.

\bibitem{Fr} M.~Freedman,
Personal communication, June 23, 1998.

\bibitem{GMM} H.~Gluck, D.~Mackenzie, and F.~Morgan,
{\em Volume-minimizing cycles in Grassmann manifolds}, Duke
Math. J. \textbf{79} (1995), 335--404.

\bibitem{G1} M.~Gromov,
{\em Filling Riemannian manifolds}, J. Differential Geom.
\textbf{18} (1983), 1--147.

\bibitem{G2} \bysame,
{\em Pseudoholomorphic curves in symplectic manifolds},
Invent. Math. \textbf{82} (1985), 307--347.

\bibitem{G3} \bysame,
{\em Systoles and intersystolic inequalities},
In: Actes de la Table Ronde de G\'{e}om\'{e}trie Diff\'{e}rentielle 
(Luminy, 1992), S\'{e}min. Congr., vol.~1, Soc. Math. France, Paris, 
1996, pp.~291--362.

\bibitem{G4} \bysame,
{\em Metric structures for Riemannian and non-Riemannian spaces},
Progr. Math., vol.~152, Birkh\"{a}user, Boston, MA, 1998.

\bibitem{G5} \bysame,
Personal communication, August 7, 1998.

\bibitem{GLP} M.~Gromov, J.~Lafontaine, P.~Pansu,
{\em Structures M\'etriques pour les Vari\'et\'es Riemanniennes},
Cedic, 1981.

\bibitem{Hil} P.~J.~Hilton,
{\em On the homotopy groups of the union of spheres},
J. London Math. Soc. \textbf{30} (1955), 154--172.

\bibitem{K1} M.~Katz,
{\em Counterexamples to isosystolic inequalities},
Geom. Dedicata \textbf{57} (1995), 195--206.

\bibitem{K2} \bysame,
{\em Systolically free manifolds}, Appendix~D in \cite{G4}.

\bibitem{LW} A.~Lundell and S.~Weingram,
{\em  The topology of CW complexes}, University Series in Higher Math.,
Van Nostrand Reinhold, New York, 1969.

\bibitem{Mi} J.~Milnor,
{\em Lectures on the $h$-cobordism theorem},
Princeton Univ. Press, Princeton, NJ, 1965.

\bibitem{P} C.~Pittet,
{\em Systoles on $S^1\times S^n$},
Differential Geom.\ Appl. \textbf{7} (1997), 139--142.

\bibitem{Se} J.-P.~Serre,
{\em  Groupes d'homotopie et classes de groupes abeliens},
Ann. of Math. \textbf{58} (1953), 258--294.

\bibitem{Su} A.~Suciu,
{\em Iterated spinning and homology spheres},
Trans. Amer. Math. Soc. \textbf{321} (1990), 145--157.

\bibitem{Sv} D.~Sullivan,
{\em Genetics of homotopy theory and the Adams conjecture},
Ann. of Math. \textbf{100} (1974), 1--79.

\bibitem{Wh1} G.~W.~Whitehead,
{\em A generalization of the Hopf invariant},
Ann. of Math. \textbf{51} (1950), 192--237.

\bibitem{Wh} \bysame,
{\em  Elements of homotopy theory}, Grad. Texts
in Math., vol.~61, Springer-Verlag, New York, 1978.

\end{thebibliography}

\end{document}